\documentclass[]{sn/interact}

\usepackage{epstopdf}% To incorporate .eps illustrations using PDFLaTeX, etc.
\usepackage[caption=false]{subfig}% Support for small, `sub' figures and tables

\usepackage[numbers,sort&compress]{natbib}% Citation support using natbib.sty
\bibpunct[, ]{[}{]}{,}{n}{,}{,}% Citation support using natbib.sty
\renewcommand\bibfont{\fontsize{10}{12}\selectfont}% Bibliography support using natbib.sty

\theoremstyle{plain}% Theorem-like structures provided by amsthm.sty
\newtheorem{theorem}{Theorem}[section]
\newtheorem{lemma}[theorem]{Lemma}

\newtheorem{proposition}[theorem]{Proposition}

\theoremstyle{definition}
\newtheorem{definition}[theorem]{Definition}

\theoremstyle{remark}
\newtheorem{remark}{Remark}[section]

\usepackage{hyperref}
\usepackage{graphicx}
\usepackage{dsfont}
\usepackage{mathtools}
\usepackage{wrapfig}
\usepackage{amssymb, amsmath, latexsym}
\usepackage{nicefrac}
\usepackage{hhline}
\usepackage{multirow}
\usepackage{comment}
\usepackage{pifont}% http://ctan.org/pkg/pifont

\usepackage[colorinlistoftodos,bordercolor=orange,backgroundcolor=orange!20,linecolor=orange,textsize=scriptsize]{todonotes}
\usepackage{algorithm}\usepackage{algpseudocode}
\usepackage[symbol]{footmisc}

\usepackage{tikz}
\usetikzlibrary{matrix, positioning, fit}

\def \R {\mathbb R}
\def \C {\mathbb{C}}
\def \perm {\mathcal{D}}
\def \wC {\widehat{\mathbb{C}}}
\def \lens {\mathcal{E}}
\def \st {\; | \;}
\def \P {\mathcal{P}_N}

\def\wtgg{\widetilde{\nabla}}

\def\O{\mathcal{O}}

\DeclareMathOperator{\spec}{spec}
\DeclareMathOperator{\diag}{diag}
\DeclareMathOperator{\conv}{conv}
\DeclareMathOperator{\Ree}{Re}
\DeclareMathOperator{\Imm}{Im}
\DeclareMathOperator{\Arg}{Arg}
\DeclareMathOperator{\Log}{Log}
\DeclareMathOperator{\supp}{supp}

\newcommand{\dA}{\Delta A}
\newcommand{\wtA}{\widetilde{A}}
\newcommand{\wtlambda}{\widetilde{\lambda}}
\newcommand{\re}{\mathcal{RE}_{\alpha}(A)}
\newcommand{\Bc}{B^c}

\newcommand{\proofref}[1]{\hyperref[proof:#1]{\ensuremath{\downarrow}}}

\begin{document}

% \articletype{ARTICLE TEMPLATE}% Specify the article type or omit as appropriate

\title{Applying Potential Theory to Lower Bounds for Gradient Methods with a Relative Error Oracle\thanks{}}

\author{
\name{
Artem Vasin \textsuperscript{a} and
Alexander Gasnikov \textsuperscript{b, a, c, d}
}
\affil{
\textsuperscript{a}Moscow Independent Research Institute of Artificial Intelligence, Moscow, Russia.
\textsuperscript{b}Innopolis University, Innopolis, Russia.
\textsuperscript{c}Steklov Mathematical Institute of the Russian Academy of Sciences, Moscow, Russia.
\textsuperscript{d}Institute for Information Transmission Problems, Moscow, Russia.
}
}

\date{Received: date / Accepted: date}
\maketitle              % typeset the header of the contribution
\begin{abstract}
    We establish a lower bound for first-order methods with relative gradient errors in smooth strongly convex optimization. The bound applies to methods satisfying the standard linear-span condition and is realized by quadratic objectives with fixed affine gradient perturbations. Our analysis combines matrix perturbation theory with potential theory. The result quantifies the joint dependence on the condition number $L/\mu$ and the relative error level $\alpha$. In particular, it identifies $\alpha=\Omega\left( 1/ \ln(L / \mu) \right)$ as a scale at which relative errors cause a non-vanishing deterioration of the accelerated convergence factor.
\end{abstract}

% \clearpage

\section{Introduction}\label{sec:Introduction}

We consider the unconstrained minimization problem:
\begin{equation}\label{eq:optim}
    \min_{x\in\mathbb{R}^n} f(x),
\end{equation}
where $f:\mathbb{R}^n\to\mathbb{R}$ is differentiable. We assume that the  minimum is attained and write $f^*=\min\limits_{x\in\mathbb{R}^n}f(x)$. Given an initial point $x^0$, let $x^*$ be a minimizer closest to $x^0$ in the Euclidean norm, and set $R=\|x^0-x^*\|_2$. Throughout, $\langle\cdot,\cdot\rangle$ and $\|\cdot\|_2$ denote the Euclidean inner product and norm, respectively.

We express $L$-smoothness, with $L>0$, through the quadratic upper model
\begin{equation}\label{eq:smooth_cond}
    (\forall x,y \in\mathbb{R}^n) \qquad f(y) \leqslant f(x) + \langle \nabla f(x), y - x \rangle
    + \frac{L}{2} \|x - y\|_2^2.
\end{equation}
We additionally assume that $f$ is $\mu$-strongly convex, with $\mu>0$:
\begin{equation}\label{eq:strong_conv}
    (\forall x,y \in\mathbb{R}^n)\qquad f(x) + \langle \nabla f(x), y - x \rangle +\frac{\mu}{2} \|x - y\|_2^2 \leqslant f(y).
\end{equation}
As a special case, we will consider a quadratic function:
\begin{equation} \label{eq:quadratic_problem}
    f_A(x) = \frac{1}{2} (x - x^*)^T A (x - x^*), \quad \mu I \leqslant A \leqslant L I, \quad A = A^T \in \R^{n \times n}.
\end{equation}

The exact gradient is assumed to be unavailable. Instead, an inexact oracle returns $\widetilde{\nabla}f(x)$ satisfying the relative error condition
\begin{equation}\label{eq:relative_noise}
    (\forall x\in\mathbb{R}^n)\qquad \|\wtgg f(x) - \nabla f(x)\|_2 \leqslant \alpha\|\nabla f(x)\|_2, \quad 0 \leqslant \alpha < 1.
\end{equation}
A special case is an oracle with a fixed matrix $\dA$:
\begin{equation} \label{eq:affine_error}
    (\forall x \in \mathbb{R}^n) \qquad \| \dA (x - x^*) \|_2 \leqslant \alpha \|\nabla f(x) \|_2, \quad \wtgg f(x) = \nabla f(x) + \dA (x - x^*).
\end{equation}
For the quadratic objective~\eqref{eq:quadratic_problem}, the oracle~\eqref{eq:affine_error} takes the form:
\begin{equation} \label{eq:affine_quadratic_model}
    \wtgg f(x) = \wtA (x - x^*), \; \wtA = A + \dA, \; \dA \in \re, \; \forall x \|\dA x \|_2 \leqslant \alpha \| A x \|_2.
\end{equation}

We consider gradient-based first-order methods, which are described in detail in~\cite{gasnikov2017universal,bubeck2015convex,beck2017first,nemirovski1995information,nesterov2018lectures,polyak1987introduction}:
\begin{equation}\label{eq:linear_first_order}
    \begin{gathered}
        x^{N + 1} \in x^0 + \operatorname{span} \left\{ \wtgg f(x^0),\wtgg f(x^1),\ldots,\wtgg f(x^{N}) \right\}, \\
        \wtgg f(x^k) \text{ is an approximation for } \nabla f(x^k).
    \end{gathered}
\end{equation}

The present paper studies lower complexity bounds under relative gradient error~\eqref{eq:relative_noise}. We study the class of first-order gradient methods~\eqref{eq:linear_first_order} and construct strongly convex quadratic instances with fixed affine gradient perturbations~\eqref{eq:affine_quadratic_model} in Theorem~\ref{thm:quadratic_lower_bound_single_point}.

\noindent \textbf{Contribution.} Theorem~\ref{thm:quadratic_lower_bound_single_point} establishes a lower bound applicable to all first-order gradient methods~\eqref{eq:linear_first_order}. For every iteration budget $N$, we construct an admissible quadratic instance in dimension $2(N + 1)$ and an initial point for which every method satisfying~\eqref{eq:linear_first_order} obeys the following:
\begin{equation*}
    \begin{gathered}
        \|x^N - x^* \|_2 \geqslant \|x^0 - x^* \|_2 \left( \frac{\widehat{L}^{\gamma} -\widehat{\mu}^{\gamma}}{\widehat{L}^{\gamma} + \widehat{\mu}^{\gamma}} \right)^{N}, \quad \gamma = \frac{\pi}{2 ( \pi - \arcsin \alpha)}, \\
        0 < \alpha < 1, \quad \widehat{L} = (1 + \alpha)L, \quad \widehat{\mu} = (1 - \alpha) \mu.
    \end{gathered}
\end{equation*}
From the perspective of the proof, our contribution is the identification of the spectral region admissible under relative gradient errors, the construction of an inscribed set with an explicit Green's function, and the application of the finite extremal construction to linear-span methods with adaptive coefficients.

% \noindent \textbf{Organization.}

\noindent \textbf{Notation.} Notation is collected in Appendix~\ref{appendix:notation}.

\section{Related work} \label{sec:related_work}

\subsection{Matrix perturbation} \label{subsec:related_matrix_perturbation}

Relative matrix perturbation theory provides a natural starting point for studying gradient errors on quadratic objectives. A relative version of the Bauer--Fike bound, established by~\cite{eisenstat1998three}, controls eigenvalue displacements relative to the eigenvalues of the unperturbed matrix. We recall the formulation used below.

\begin{theorem}[Corollary of Bauer-Fike Theorem, Corollary 2.2~\cite{eisenstat1998three}] \label{thm:bauer-fike}
    Let $A\in\C^{n\times n}$ be diagonalizable and non-singular, with
    \begin{equation*}
        A=X\Lambda X^{-1},\qquad \Lambda=\diag(\lambda_1,\ldots,\lambda_n).
    \end{equation*}
    For every $E\in\C^{n\times n}$ and every
    $z\in\spec(A+E)$,
    \begin{equation*}
        \min_{1\leq j\leq n} \frac{|z - \lambda_j|}{|\lambda_j|}
        \leq \kappa_2(X) \| A^{-1} E \|_2, \quad \kappa_2(X) = \| X \|_2 \| X^{-1} \|_2.
    \end{equation*}
\end{theorem}

Theorem~\ref{thm:bauer-fike} connects the relative gradient error model~\eqref{eq:affine_quadratic_model} to the geometry of admissible complex spectra.

\subsection{Potential theory} \label{subsec:related_potential_theory}

Potential theory supplies the tools needed to translate this spectral geometry into polynomial approximation bounds. We use the framework developed in~\cite{ransford1995potential}, particularly Green's functions, their conformal invariance, and the Bernstein--Walsh inequality.

\begin{theorem}[Bernstein's Lemma, Theorem 5.5.7~\cite{ransford1995potential}] \label{thm:bernstein_theorem}
    Let $K$ be a non-polar (Definition~\ref{def:polar_set}) compact subset of $\C$, and let $D$ be the component of $\wC \setminus K$ containing $\infty$. If $q$ is a polynomial of degree $n \geqslant 1$, then
    \begin{equation*}    
        (\forall z \in D \setminus \{\infty\}) \quad \left( \frac{|q(z)|}{\|q\|_K} \right)^{1/n} \leqslant e^{g_D(z, \infty)},
    \end{equation*}
    where $g_D$ denotes the Green's function of $D$.
\end{theorem}

%Since a first-order method for a quadratic problem with affine error~\eqref{eq:affine_quadratic_model} can be represented as the application of an operator $q(\wtA)$ to the initial residual $x^0 - x^*$ (Lemma~\ref{lem:polynomial_expression}), the proposed theorem can be used for the case where $K$ is a set containing the spectrum of the matrix $\wtA$. Details of potential theory are presented in Appendix~\ref{appendix:complex_potential}. 

Potential theory has a long history in the analysis of iterative methods. The works of~\cite{nevanlinna2012convergence, driscoll1998potential} describe how extremal polynomials, Green's functions, and conformal mappings determine convergence rates for Krylov subspace and semi-iterative methods. Related conformal mapping techniques also appear in GMRES convergence estimates~\cite{beckermann2005some}. Similar ideas have also been used by~\cite{azizian2020accelerating} to derive spectral-shape lower bounds for oblivious first-order methods in smooth games. The authors of~\cite{kelner2022big} use potential theory to establish lower bounds for multiscale optimization. In~\cite{shugart2025min} Green's functions and conformal mappings are used to characterize optimal convergence rates for quadratic min-max problems and related variational inequalities.
%The minimax principle underlying Theorem~\ref{thm:existance_points_z_k} is closely related to the equality between worst-case and ideal approximation for normal matrices established by~\cite{greenbaum1994max, joubert1994robust}.

\subsection{Gradient methods} \label{subsec:related_gradient_methods}

Classical convergence guarantees for methods with exact gradients are presented in~\cite{nesterov2018lectures}. For an $L$-smooth~\eqref{eq:smooth_cond} convex objective, gradient descent and accelerated gradient methods achieve objective residuals of order $\frac{LR^2}{N}$ and $\frac{LR^2}{N^2}$, respectively. Under $\mu$-strong convexity~\eqref{eq:strong_conv}, their standard iteration bounds become
\begin{equation*}
    \O \left( \frac{L}{\mu} \log\left(\frac{LR^2}{\varepsilon}\right) \right) \quad \text{and} \quad \O \left( \sqrt{\frac{L}{\mu}} \log\left(\frac{LR^2}{\varepsilon}\right) \right),
\end{equation*}
for the target $f(x^N)-f(x^*)\le\varepsilon$. Classical lower complexity bounds establish the optimal worst-case orders associated with acceleration in the standard first-order
oracle setting~\eqref{eq:linear_first_order}. The convergence rate of gradient descent under relative gradient errors was investigated by~\cite{de2017worst, de2020worst, vernimmen2025worst} using the PEP technique~\cite{goujaud2024pepit}.

For inexact gradients with relative error~\eqref{eq:relative_noise}, the relevant question is how these guarantees change under the relative error condition. The first positive results of convergence of gradient descent with such an imprecise oracle were given in~\cite{polyak1987introduction}.

A frequency-domain analysis of inexact gradient methods using stability criteria from control theory was developed in~\cite{gannot2022frequency}. In that paper, an accelerated method is introduced, and its convergence is established for strongly convex~\eqref{eq:strong_conv} and smooth~\eqref{eq:smooth_cond} functions, for $\alpha = O\left(\sqrt{\frac{\mu}{L}}\right)$:
\begin{equation*}
    \|x^N - x^*\|_2 = \O \left(\left(1 - C \sqrt{\frac{\mu}{L}}\right)^N \|x^0 - x^*\|_2 \right), \quad 0 < C < 1 - \text{ an absolute constant}.
\end{equation*}
The same result was obtained using a restarting technique by~\cite{kornilov2025intermediate}.

The authors of~\cite{vasin2025solving} proposed a single-loop relative-error intermediate gradient method, RE-AGM. Their results imply:
\begin{equation*}
    \O \left( \left(\sqrt{\frac{L}{\mu}} + \alpha \frac{L}{\mu} \right)\log\left(\frac{LR^2}{\varepsilon}\right) \right),
\end{equation*}
oracle complexity for $0\leqslant \alpha \leqslant 1/3$. Consequently, acceleration is preserved when $\alpha = \O((L / \mu)^{-1/2})$, while larger admissible errors yield intermediate linear rates, reaching the gradient descent order for constant $\alpha$. The paper also extends relative-error analysis to min--max problems and operator methods.

More recently, the intermediate method proposed in~\cite{vasin2026adaptive} improved the dependence on the relative error. Its complexity is:
\begin{equation*}
    \O \left( \left(\sqrt{\frac{L}{\mu}} + \alpha^2 \frac{L}{\mu} \right)\log\left(\frac{LR^2}{\varepsilon}\right) \right), \quad 0 \leqslant \alpha \leqslant \frac{1}{3}.
\end{equation*}
The accelerated order therefore persists up to $\alpha = \O((L / \mu)^{-1/4})$. The authors of~\cite{tyurin2026longtails} obtained related oracle-complexity guarantees within the \textit{LongTails} framework for formally verified optimization results. Their strongly convex guarantee covers every $0 \leqslant \alpha < 1$ and has complexity:
\begin{equation*}
    \O \left( \left(\sqrt{\frac{L}{(1 - \alpha)\mu}} + \alpha^2 \frac{L}{(1 - \alpha) \mu} \right)\log\left(\frac{LR^2}{\varepsilon}\right) \right),
\end{equation*}
For smooth convex objectives, they also established the following oracle-complexity bound for obtaining a last iterate satisfying $f(x^N) - f^* \leqslant \varepsilon$:
\begin{equation*}
    \O \left( \sqrt{\frac{L R^2}{(1 - \alpha)\varepsilon}} + \alpha^2 \frac{L R^2}{(1 - \alpha) \varepsilon} \right).
\end{equation*}

\begin{remark} \label{rem:error_in_previous_lb}
    The proof of Theorem 10.1 in the first version of the preprint~\cite{vasin2026lower} contained an incorrect reduction involving the inexact oracle model provided in~\cite{devolder2013first}. The proposed $(\delta, L, \mu)$-oracle model considered in~\cite{devolder2013first} is more general than the model of absolute gradient error $\| \wtgg f(x) - \nabla f(x) \|_2 \leqslant \delta, \; \delta > 0$. In preprint~\cite{vasin2026lower}, a lower bound for the $(\delta, L, \mu)$ model was applied to the absolute gradient error model, and subsequently, the absolute error was converted into a relative error using a stopping criterion:
    \begin{equation*}
        \begin{gathered}
            \|\wtgg f(x) \|_2 \geqslant K \delta \Longrightarrow \frac{\|\nabla f(x) - \wtgg f(x) \|_2}{\| \nabla f(x) \|_2} \leqslant \frac{\delta}{(K - 1) \delta} = \frac{1}{K - 1}, \quad K > 2, \\
            \text{else } \|\wtgg f(x) \|_2 < K \delta \Longrightarrow \text{ stop method}. 
        \end{gathered}
    \end{equation*}
    Consequently, a lower bound for the former does not automatically apply to the latter: one must additionally show that the hard instances satisfy the prescribed absolute gradient-error bound. Thus, the cited argument does not establish the claimed lower bound. The optimal dependence of the worst-case complexity on the relative error level and the condition number remains unresolved.
\end{remark}

\section{Ideas behind the result} \label{sec:ideas_behind}

\begin{figure}[htbp]
\centering
\begin{tikzpicture}[x=1cm,y=1cm,font=\small]
  % Illustration parameters.
  \pgfmathsetmacro{\figmu}{3}
  \pgfmathsetmacro{\figL}{6}
  \pgfmathsetmacro{\figalpha}{0.45}

  % TikZ trigonometric functions use degrees.
  \pgfmathsetmacro{\figtheta}{asin(\figalpha)}
  \pgfmathsetmacro{\figa}{(1-\figalpha)*\figmu}
  \pgfmathsetmacro{\figb}{(1+\figalpha)*\figL}
  \pgfmathsetmacro{\figr}{(\figb-\figa)/(2*\figalpha)}
  \pgfmathsetmacro{\figs}{sqrt(1-\figalpha*\figalpha)}

  % Tangency points of the endpoint disks.
  \pgfmathsetmacro{\figxm}{\figmu*(1-\figalpha*\figalpha)}
  \pgfmathsetmacro{\figym}{\figalpha*\figmu*\figs}
  \pgfmathsetmacro{\figxL}{\figL*(1-\figalpha*\figalpha)}
  \pgfmathsetmacro{\figyL}{\figalpha*\figL*\figs}

  % Exact boundary of the continuous union of disks:
  % two common tangents and arcs of the endpoint disks.
  \path[fill=blue!9,draw=blue!65!black,line width=0.9pt]
    (\figxm,\figym) -- (\figxL,\figyL)
    arc[start angle={90+\figtheta},
        end angle={-90-\figtheta},
        radius={\figalpha*\figL}]
    -- (\figxm,-\figym)
    arc[start angle={-90-\figtheta},
        end angle={-270+\figtheta},
        radius={\figalpha*\figmu}]
    -- cycle;

  % Representative disks from the continuous union.
  \foreach \lam in {3,4,5,6}{
    \draw[blue!35,densely dashed,line width=0.4pt]
      (\lam,0) circle[radius={\figalpha*\lam}];
  }

  % Exact circular lens E(a,b,theta).
  \path[fill=orange!30,draw=orange!70!black,line width=1pt]
    (\figa,0)
    arc[start angle={90+\figtheta},
        end angle={90-\figtheta},
        radius=\figr]
    arc[start angle={-90+\figtheta},
        end angle={-90-\figtheta},
        radius=\figr]
    -- cycle;

  % Coordinate axes.
  \draw[->,black!60] (-0.15,0) -- (9.65,0)
    node[right] {$\Ree z$};
  \draw[->,black!60] (0,-3.05) -- (0,3.15)
    node[above] {$\Imm z$};
  \node[below left] at (0,0) {$0$};

  % The original real spectral interval.
  \draw[black,line width=1.8pt]
    (\figmu,0) -- (\figL,0);

  \foreach \xx in {\figa,\figmu,\figL,\figb}{
    \fill (\xx,0) circle[radius=1.5pt];
  }

  \node[above left=3pt] at (\figa,0) {$(1 - \alpha) \mu$};
  \node[below=4pt] at (\figmu,0) {$\mu$};
  \node[below=4pt] at (\figL,0) {$L$};
  \node[above right=3pt] at (\figb,0) {$(1 + \alpha) L$};

  \node[text=orange!60!black] at (4.8,0.52)
    {$\lens$};
  \node[text=blue!65!black] at (6.2,2.25) {$D$};
\end{tikzpicture}
\caption{$D$ is the set of admissible eigenvalues for the perturbed matrix $\wtA$, and $\lens$ is the set from which eigenvalues are selected for $\wtA$.}
\label{fig:spectral-lens}
\end{figure}
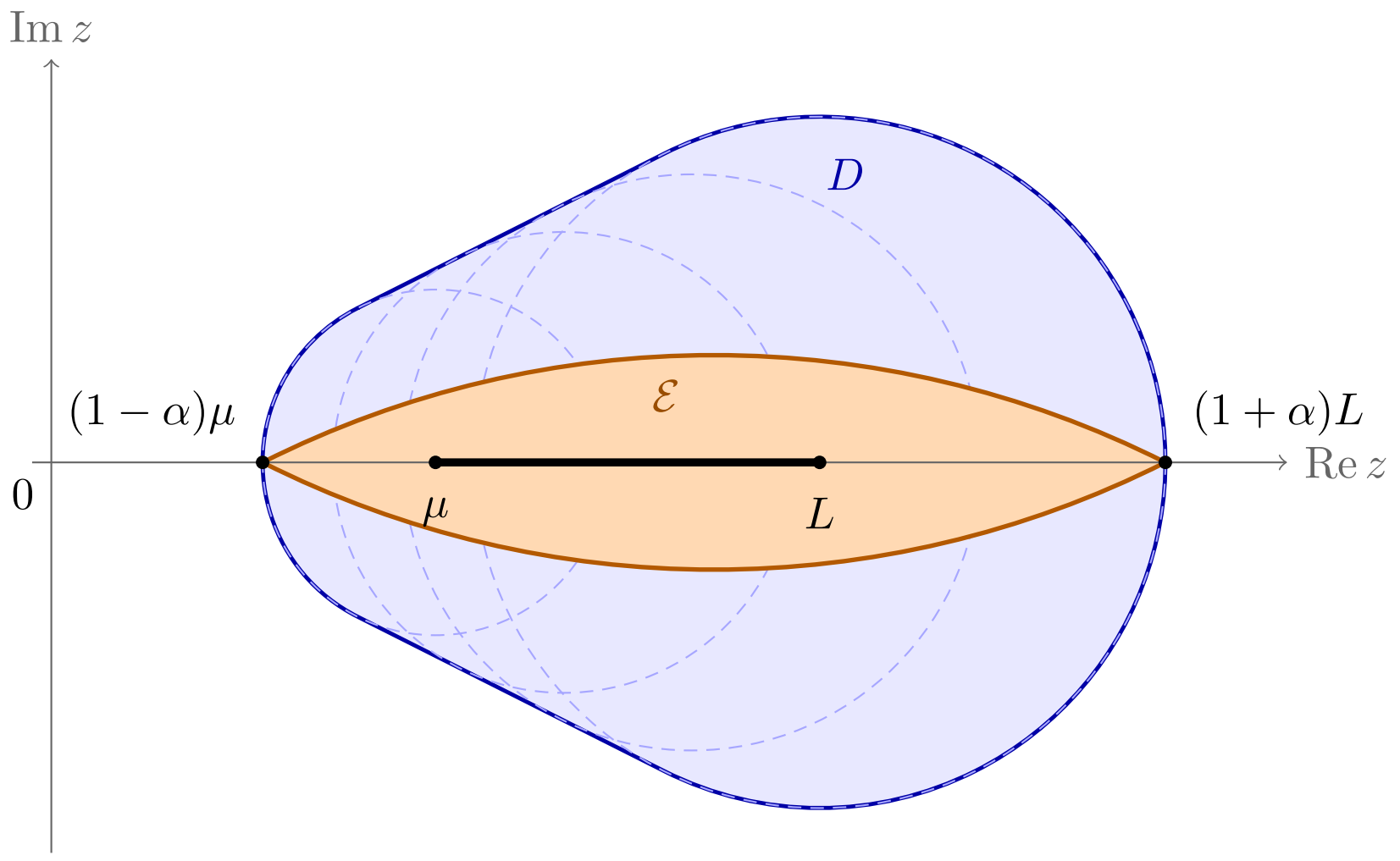

The proof of Theorem~\ref{thm:quadratic_lower_bound_single_point} reduces the convergence of first-order methods to polynomial approximation on a suitable complex set. We consider the quadratic function and affine oracle~\eqref{eq:affine_quadratic_model}. Lemma~\ref{lem:polynomial_expression} shows the iterates of any first-order method satisfying~\eqref{eq:linear_first_order} admit the representation:
\begin{equation*}    
    x^N - x^* = P_N(\wtA) (x^0 - x^*), \quad
    P_N \in \P = \{ P \in \R[z]: \deg P \leqslant N, \; P(0) = 1 \}.
\end{equation*}
Thus, a lower bound can be obtained by constructing an admissible matrix and an initial vector for which every such polynomial leaves a substantial value of $\|x^N - x^* \|_2$. The consequence (Lemma~\ref{lem:admissible_errors}) of the Bauer--Fike bound (Theorem~\ref{thm:bauer-fike}) gives the spectral enclosure:
\begin{equation*}
    \spec(\wtA) \subseteq D := \bigcup_{\lambda\in[\mu,L]} \{z \in \C: |z - \lambda| \leqslant \alpha \lambda \}.
\end{equation*}
We choose a lens $\lens \subseteq D$. The definition and properties of $\lens$ can be found in Appendix~\ref{appendix:lens_prop}. In particular, its boundary consists of two circular arcs joining $(1 - \alpha) \mu$ and $(1 + \alpha) L$. Figure~\ref{fig:spectral-lens} illustrates this inclusion. Appendix~\ref{appendix:lens_prop} also provides the conformal map:
\begin{equation*}
     \Phi(z) = \frac{\left(\frac{a - z}{b - z} \right)^{\gamma} + 1}{\left(\frac{a - z}{b - z} \right)^{\gamma} - 1}, \qquad \gamma = \frac{\pi}{2 ( \pi - \arcsin \alpha)}, \quad a = (1 - \alpha) \mu, \quad b = (1 + \alpha) L,
\end{equation*}
from $\wC \setminus \lens$ onto the exterior of the unit disk $ \Bc = \{ z \in \C \st |z| > 1 \} \cup \{ \infty \}$, with $\Phi(\infty)=\infty$. The Lemma~\ref{lem:green_function_using_conf_map} shows that $\ln |\Phi(z) |$ will be the Green's function $g(z, \infty)$ for $\wC \setminus \lens$.

Using Danskin's Theorem (Theorem~\ref{thm:danskin_sub}) and Caratheodory's Theorem (Theorem~\ref{thm:caratheodory_thm}), we select in Theorem~\ref{thm:existance_points_z_k} points $z_1,\ldots,z_{N+1}\in \lens$ and weights $\omega_j\ge0$, with $\sum\limits_{j = 1}^{N + 1}\omega_j=1$, such that:
\begin{equation*}
    (\forall P \in \P) \quad \sum_{j=1}^{N+1} \omega_j|P(z_j)|^2 \geqslant M^2, \quad
    M = \inf_{P \in \P}\max_{z\in \lens} |P(z)|.
\end{equation*}
For each node, the inclusion $\lens \subseteq D$ supplies $\lambda_j \in [\mu,L]$ with $|z_j-\lambda_j | \leqslant \alpha \lambda_j$. Representing multiplication by $z_j$ as a real $2 \times 2$ block produces $\widetilde A$ with eigenvalues $z_j, \overline{z_j}$, while the corresponding blocks $\lambda_j I_2$ define $A$. With a suitable choice of the initial point, one can obtain:
\begin{equation*}
    \|x^N - x^*\|_2^2 =R^2 \sum_{j=1}^{N+1} \omega_j|P_N(z_j)|^2 \geqslant R^2 M^2.
\end{equation*}
Applying Bernstein's lemma (Theorem~\ref{thm:bernstein_theorem}) at $z=0$ gives:
\begin{equation*}
    M \geqslant | \Phi(0) |^{-N} = \left( \frac{b^\gamma-a^\gamma}{b^\gamma+a^\gamma} \right)^N.
\end{equation*}
Combining the last two bounds proves Theorem~\ref{thm:quadratic_lower_bound_single_point}. Details can be found in Theorem~\ref{thm:quadratic_lower_bound_single_point} and Appendix~\ref{appendix:lower_bound}.

\section{Lower bound} \label{sec:lower_bound}

Finally, we move on to the proof of the main theorem and its discussion. We assume that $\alpha > 0$, since otherwise there is no noise (see Remark~\ref{rem:lower_bound_alpha=0}).

\begin{theorem} \label{thm:quadratic_lower_bound_single_point}
    Let $0 < \alpha < 1, \; 0 < \mu < L$, $N \geqslant 1$ and $R > 0$. There exist $x^0 \in \R^{2(N + 1)}$, $0 < \mu I \leqslant A \leqslant L I$ and $\dA \in \re$~\eqref{eq:affine_quadratic_model} ($\forall x \; \|\dA x \|_2 \leqslant \alpha \| A x \|_2$), such that for all first order algorithms~\eqref{eq:linear_first_order} $\mathcal{A}$ producing $x^N$ (the point at the $N$-th iteration) the following lower bound holds:
    \begin{equation*}
        \|x^N - x^* \|_2 \geqslant R \left( \frac{\widehat{L}^{\gamma} -\widehat{\mu}^{\gamma}}{\widehat{L}^{\gamma} + \widehat{\mu}^{\gamma}} \right)^{N}, \quad \gamma = \frac{\pi}{2 ( \pi - \arcsin \alpha)}, \quad \widehat{L} = (1 + \alpha)L, \; \widehat{\mu} = (1 - \alpha) \mu.
    \end{equation*}
    where $f(x) = \frac{1}{2} x^T A x, \; f^* = 0, \; x^* = 0, \; \|x^0 - x^* \|_2 = R$.
\end{theorem}

% \left(1 - \left(\displaystyle \frac{1 - \alpha}{1 + \alpha} \cdot \frac{\mu}{L}\right)^{\gamma}\right)^N

\begin{remark} \label{rem:lower_bound_log}
    Let us consider $0 < \mu < \frac{1}{4} L$ and $\alpha \geqslant \displaystyle \frac{2 \pi \ln c_1}{\ln \nicefrac{L}{\mu}}, \; c_1 > 1$. Let us denote $t = \frac{1 - \alpha}{1 + \alpha} \cdot \frac{\mu}{L}$. Then applying Theorem~\ref{thm:quadratic_lower_bound_single_point}:
    \begin{equation*}
         \|x^N - x^* \|_2 \geqslant R \left( \frac{\widehat{L}^{\gamma} -\widehat{\mu}^{\gamma}}{\widehat{L}^{\gamma} + \widehat{\mu}^{\gamma}} \right)^{N}
         = R \left( \frac{1 - t^{\gamma}}{1 + t^{\gamma}} \right)^{N} \geqslant R \left(1 - 2 t^{\gamma} \right)^{N}.
    \end{equation*}
    Using:
    \begin{eqnarray*}
        \gamma & = & \frac{\pi}{2 ( \pi - \arcsin \alpha)} \overset{\arcsin \alpha \geqslant \alpha}{\geqslant} \frac{1}{2(1 - \nicefrac{\alpha}{\pi})} \geqslant \frac{1}{2} + \frac{\alpha}{2 \pi} \geqslant \frac{1}{2} + \frac{\ln c_1}{\ln \nicefrac{L}{\mu}}
        \\
        & \geqslant & \frac{1}{2} -\frac{\ln c_1}{\ln t} \geqslant \frac{1}{2} + \log_t{\frac{1}{c_1}}
    \end{eqnarray*}
    Then:
    \begin{equation*}
        \displaystyle 1 - 2 t^{\gamma} \geqslant 1 - 2 \sqrt{t} \cdot t^{\displaystyle \log_t{\frac{1}{c_1}}} = 1 - \frac{2}{c_1} \sqrt{t}.
    \end{equation*}
    That is:
    \begin{equation*}
         \|x^N - x^* \|_2 \geqslant R \left(1 - \displaystyle \frac{2}{c_1} \sqrt{\frac{1 - \alpha}{1 + \alpha} \cdot \frac{\mu}{L}} \right)^N.
    \end{equation*}
    Thus, Theorem~\ref{thm:quadratic_lower_bound_single_point} shows that relative errors of size $\alpha = \Omega \left( \displaystyle \frac{1}{\ln \nicefrac{L}{\mu}} \right)$ already cause a non-vanishing deterioration relative to the classical exact-gradient convergence factor. We can also consider the case $\alpha(\kappa)$, where $\kappa = \frac{L}{\mu}$. When $\alpha(\kappa)\ln\kappa\to+\infty$ as $\kappa\to+\infty$, Theorem~\ref{thm:quadratic_lower_bound_single_point} rules out a uniform worst-case convergence guarantee with the classical accelerated $\sqrt{\kappa}$ dependence for methods satisfying~\eqref{eq:linear_first_order}.
\end{remark}

\begin{remark} \label{rem:convex_lb}
Theorem~\ref{thm:quadratic_lower_bound_single_point} also yields a lower bound for smooth convex optimization, where no positive strong convexity parameter is prescribed. Choose $N \geqslant 2$ and $\mu$ satisfying:
\begin{equation*}    
    N > \left(\frac{ 1 +\alpha}{1 - \alpha}\right)^{\gamma},
    \quad
    \gamma = \frac{\pi}{2 (\pi - \arcsin\alpha)},
    \quad
    \mu = \frac{1 + \alpha}{1 - \alpha} \frac{L}{N^{1/\gamma}}.
\end{equation*}
The condition on $N$ ensures that $0 < \mu < L$, so Theorem~\ref{thm:quadratic_lower_bound_single_point} applies. The resulting quadratic objective is, in particular, $L$-smooth and convex. Moreover, $\left(\frac{(1 - \alpha)\mu}{(1 + \alpha) L}\right)^\gamma = \frac{1}{N}$. Thus, there exist such an objective in dimension $2(N + 1)$, a fixed affine oracle satisfying~\eqref{eq:affine_error}, and an initial point with $\|x^0 - x^*\|_2 = R > 0$ for which every method satisfying~\eqref{eq:linear_first_order} obeys:
\begin{equation*}
    f(x^N)-f^* \geqslant \frac{\mu}{2}\|x^N - x^*\|_2^2
    \overset{\text{Theorem~\ref{thm:quadratic_lower_bound_single_point}}}{\geqslant} \frac{\mu R^2}{2} \left(\frac{N - 1}{N + 1}\right)^{2N} \geqslant \frac{1}{162} \frac{1 + \alpha}{1 - \alpha} \frac{LR^2}{N^{2 (1 - \frac{\arcsin \alpha}{\pi})}}.
\end{equation*}
Here we used $\left(\frac{N - 1}{N + 1}\right)^{2N} \geqslant \frac{1}{81}$ for $N \geqslant 2$
and $\frac{1}{\gamma} = 2 \left(1 - \frac{\arcsin\alpha}{\pi} \right)$. Consequently, for fixed $\alpha$ and sufficiently small $\varepsilon>0$, the worst-case number of oracle calls required
to guarantee $f(x^N) - f^* \leqslant \varepsilon$ is at least:
\begin{equation*}
    \Omega \left( \left( \frac{1 + \alpha}{1 -\alpha} \frac{LR^2}{\varepsilon} \right)^\gamma \right), \quad \gamma = \frac{\pi}{2 (\pi - \arcsin\alpha)}.
\end{equation*}
In particular, any fixed positive relative error rules out a uniform worst-case $\O \left(\frac{LR^2}{N^2} \right)$ guarantee over this class.
\end{remark}

\begin{remark} \label{rem:lower_bound_alpha=0}
    A similar lower bound for the exact-gradient setting is given in~\cite{nesterov2018lectures} for $\mu$-strongly convex~\eqref{eq:strong_conv} and $L$-smooth~\eqref{eq:smooth_cond} functions using an exact gradient oracle:
    \begin{equation} \label{eq:nesterov_lb}
        \|x^N - x^* \|_2 \geqslant \| x^0 - x^* \|_2 \left( \frac{\sqrt{L} - \sqrt{\mu}}{\sqrt{L} + \sqrt{\mu}} \right)^N.
    \end{equation}
    Note that the estimate from Theorem~\ref{thm:quadratic_lower_bound_single_point} tends to the bound~\eqref{eq:nesterov_lb} when $\alpha \to 0+$.
\end{remark}
\section{Conclusion} \label{sec:conslusion}

We established a lower bound for first-order methods satisfying the linear-span condition under relative gradient errors~\eqref{eq:relative_noise}. By combining matrix perturbation theory with potential theory, we constructed strongly convex quadratic instances with fixed affine oracles that quantify the effect of relative errors on convergence.

A gap remains between the available upper bounds~\cite{vasin2026adaptive, tyurin2026longtails} and the lower bound established here. Writing $\kappa = \frac{L}{\mu}$, accelerated convergence is achievable when $\alpha = \O (\kappa^{-1/4})$, whereas our lower bound (Theorem~\ref{thm:quadratic_lower_bound_single_point}) rules out uniform accelerated convergence when $\alpha\ln\kappa \to +\infty$ as $\kappa \to +\infty$. A detailed investigation of the sharpness of the bound, including the class of quadratic objectives with fixed affine oracles~\eqref{eq:affine_quadratic_model}, is left for future work. This investigation will aim to clarify the optimal dependence on the relative error level and the condition number.

\section*{AI Use Statement}

The authors used GPT-6 to generate ideas, check proofs, and improve the English text (spelling and grammar). All the proofs were written by the authors themselves.

\bibliographystyle{sn/tfs}
\bibliography{literature}

\begin{appendix}

\section{Notation} \label{appendix:notation}

\begin{definition} \label{def:extended_complex_plane}
    We denote the extended complex plane as $\wC = \C \cup \{ \infty \}$. Topologically, it is homeomorphic to the two-dimensional Riemann sphere.
\end{definition}

\begin{definition} \label{def:sup_norm}
    Let $A \subset \C$ and $f: \C \to \C$:
    \begin{equation*}
        \|f\|_A = \sup_{z \in A} |f(z)|.
    \end{equation*}
\end{definition}

\begin{definition} \label{def:polynom_N}
    \begin{equation*}
        \P = \{P \in \R[z]: \deg P \leqslant N,\ P(0) = 1\}.
    \end{equation*}
\end{definition}

\begin{definition} \label{def:re}
    \begin{equation*}
        \re = \{ \dA \in \R^{n \times n} \st  \forall x \|\dA x \|_2 \leqslant \alpha \| A x \|_2 \}
    \end{equation*}
\end{definition}

\begin{definition} \label{def:subdiff}
    The subdifferential of function $ f $ at $x_0$ is denoted by $\partial f(x_0)$ and is defined as follows:
    \begin{equation}\label{eq:subdiff}
        \partial f(x_0) = \left\{ g \in \mathbb{R}^n \mid f(x) \ge f(x_0) + \langle g, x - x_0 \rangle \quad \forall x \in \mathbb{R}^n \right\}.
    \end{equation}
\end{definition}

\section{Convex analysis} \label{appendix:convex_analysis}

\begin{theorem}[Caratheodory's Theorem, Proposition 1.2.1~\cite{bertsekas2009convex}] \label{thm:caratheodory_thm}
    Let \(X\) be a non-empty subset of $\R^n$. Every vector from $\conv(X)$ can be represented as a convex combination of no more than $n + 1$ vectors from $X$.
\end{theorem}

\begin{theorem}[Theorem 3.1.20~\cite{nesterov2018lectures}] \label{thm:subdiff_optim}
    Let $X^* = \operatorname*{arg\,min}\limits_{x} f(x)$ and let the function $f$ be subdifferentiable. $x^* \in X^*$ if and only if:
    \begin{equation*}
        0 \in \partial f(x^*).
    \end{equation*}
\end{theorem}

\begin{theorem}[Danskin's Theorem, Proposition A.3.2~\cite{bertsekas2009convex}] \label{thm:danskin_sub}
    Let $Z \subset \R^m$ be a compact nonempty set, and let $\phi : \R^n \times Z \to \R$ be continuous and such that $\phi(\cdot, z) : \R^n \to \R$ is convex for each $z \in Z$. The function $f : \R^n \to \R$ given by:
    \begin{equation} \label{eq:danskin_def}
        f(x) = \max_{z \in Z} \phi(x, z)
    \end{equation}
    is convex. Let us define $Z(x)$ as the set of maximizing points in~\eqref{eq:danskin_def}
    \begin{equation*}
        Z(x) = \left\{ \overline{z} \;\middle|\; \phi(x, \overline{z}) = \max_{z \in Z} \phi(x, z) \right\}.
    \end{equation*}
    If $\phi(\cdot, z)$ is differentiable for all $z \in Z$ and $\nabla_x \phi(x, \cdot)$ is continuous on $Z$ for each $x$, then
    \begin{equation*}
        \partial f(x) = \conv \{\nabla_x \phi(x, z) \mid z \in Z(x)\}, \qquad \forall x \in \R^n.
    \end{equation*}
\end{theorem}

\section{Quadratic targets and affine error} \label{appendix:quadratic_misc}

\begin{lemma} \label{lem:admissible_errors}
    Let $0 < \mu I \leqslant A \leqslant L I$, $\dA \in \re$~\eqref{eq:affine_quadratic_model}, $\wtA = A + \dA$. Then:
    \begin{equation*}
        (\forall \wtlambda \in \spec \wtA) (\exists \lambda \in \spec A) \quad |\wtlambda - \lambda | \leqslant \alpha \lambda.
    \end{equation*}
\end{lemma}
\begin{proof}
    Using Theorem~\ref{thm:bauer-fike} for $\wtA^T = A + \dA^T$ one can obtain for $\wtlambda \in \spec \wtA^T = \spec \wtA$:
    \begin{equation*}
        \min_{\lambda \in \spec A} \frac{|\wtlambda - \lambda|}{\lambda} \leqslant \kappa_2(U) \| A^{-1} \dA^T \|_2, \quad \kappa_2(U) = \| U \|_2 \|U^{-1} \|_2.
    \end{equation*}
    Since $A = U \Lambda U^T$ is real symmetric, $U$ is orthogonal, hence $\kappa_2(U) = 1$. From the definition of $\dA$~\eqref{eq:affine_quadratic_model}:
    \begin{equation*}
        (\forall x \in \R^n) \; \|\dA x \|_2 \leqslant \alpha \|A x\|_2 \Longrightarrow (\forall x \in \R^n) \; \|\dA A^{-1} x \|_2 \leqslant \alpha \|x\|_2,
    \end{equation*}
    that is $\| A^{-1} \dA^T \|_2 = \|\dA A^{-1} \|_2 \leqslant \alpha$. Then:
    \begin{equation*}
        \min_{\lambda \in \spec A} \frac{|\wtlambda - \lambda|}{\lambda} \leqslant \alpha \Longrightarrow (\exists \lambda \in \spec A) \quad |\wtlambda - \lambda | \leqslant \alpha \lambda.
    \end{equation*}
\end{proof}

\begin{lemma} \label{lem:T_operator_isometry}
    Consider the homomorphism: $T: \C \to \R^{2 \times 2}$, such that:
    \begin{equation*}
        T(z) = \begin{pmatrix} \Ree z & - \Imm z \\ \Imm z & \Ree z \end{pmatrix}.
    \end{equation*}
    Then:
    \begin{equation*}
        (\forall v \in \R^2) \quad \| T(z) v \|_2^2 = |z|^2 \|v \|_2^2.
    \end{equation*}
\end{lemma}
\begin{proof}
    \begin{equation*}
        \|(T(z)) v \|_2^2 
        \overset{z = x + i y}{=} (x v_1 - y v_2)^2 + (y v_1 + x v_2)^2 = (x^2 + y^2) (v_1^2 + v_2^2) =  |z|^2 \|v \|_2^2.
    \end{equation*}
\end{proof}

\begin{lemma} \label{lem:lambda_and_z}
    Let $z_1, \dots z_m \in \C, \; 0 < \lambda_1, \dots \lambda_m \in \R$, such that $\forall k, \; |z_k - \lambda_k| \leqslant \alpha \lambda_k$ for $0 \leqslant \alpha < 1$ and $P \in \R[z]$ be a polynomial. Consider the homomorphism: $T: \C \to \R^{2 \times 2}$, such that:
    \begin{equation*}
        T(z) = \begin{pmatrix} \Ree z & - \Imm z \\ \Imm z & \Ree z \end{pmatrix},
    \end{equation*}
    and $A = \diag(T(\lambda_1), T(\lambda_2), \dots T(\lambda_m)), \; I_2 = \diag(1, 1), \; \wtA = \diag(T(z_1), T(z_2), \dots T(z_m))$. Then $(\forall v \in R^{2m}) \; \|(\wtA - A) v \|_2 \leqslant \alpha \|A v\|_2$, that is $(\wtA - A) \in \re$~\eqref{eq:affine_quadratic_model}.
\end{lemma}
\begin{proof}
    Since $A$ and $\wtA$ are block-diagonal, it suffices to consider one block $T(z_k) - T(\lambda_k)$ for $v = (v_1, v_2) \in \R^2$:
    \begin{eqnarray*}
        \|(T(z_k) - T(\lambda_k)) v \|_2^2 
        & = & \|T(z_k - \lambda_k) v \|_2^2 \overset{\text{Lemma}~\ref{lem:T_operator_isometry}}{=} 
        |z_k - \lambda_k|^2 \|v \|_2^2 \leqslant \alpha^2 \lambda_k^2 \| v \|_2^2.
    \end{eqnarray*}
    That is $\|(T(z_k) - T(\lambda_k)) v \|_2 \leqslant \alpha \| T(\lambda_k) v \|_2$, then combining all blocks: $(\forall v \in R^{2m}) \; \|(\wtA - A) v \|_2 \leqslant \alpha \|A v\|_2$.
\end{proof}

\begin{lemma} \label{lem:permissible_spectrum_envelope}
    Let $0 < \mu \leqslant L, 0 \leqslant \alpha < 1$ and:
    \begin{equation*}
        \perm = \{ z \in \C \st \exists \lambda \in [\mu; L], \; |z - \lambda | \leqslant \alpha \lambda \}
    \end{equation*}
    Then $\perm$ can be expressed as:
    \begin{equation*}
        \perm = \{ z \in \C \st |\Imm z | \leqslant u(\Ree z), (1 - \alpha)\mu \leqslant \Ree z \leqslant (1 + \alpha) L \},
    \end{equation*}
    \begin{equation*}
        u(x) = 
        \begin{cases}
            \sqrt{\alpha^2 \mu^2 - (x - \mu)^2}, \quad \frac{x}{1 - \alpha^2} < \mu, \\
            \frac{\alpha}{\sqrt{1 - \alpha^2}} x, \quad \mu \leqslant \frac{x}{1 - \alpha^2} \leqslant L, \\
            \sqrt{\alpha^2 L^2 - (L - x)^2}, \quad \frac{x}{1 - \alpha^2} > L.
        \end{cases}
    \end{equation*}
\end{lemma}
\begin{proof}
    The inequality $|z - \lambda| \leqslant \alpha \lambda$ for $z = x + iy$ can be written as:
    \begin{equation*}
        (x - \lambda)^2 + y^2 \leqslant \alpha^2 \lambda^2.
    \end{equation*}
    A point $z = x + iy$ belongs to $\perm$, if the inequality above holds for some $\lambda$, that is admissible points satisfy:
    \begin{equation*}
        y^2 \leqslant \sup_{\mu \leqslant \lambda \leqslant L} \left(-(x - \lambda)^2 + \alpha^2 \lambda^2 \right).
    \end{equation*}
    Consider the parabola $g(\lambda) = -(x - \lambda)^2 + \alpha^2 \lambda^2$. Its vertex can be found as $\lambda_0 = \frac{x}{1 - \alpha^2}$. Then the maximum of this function is achieved at $\lambda_0$, if $\mu \leqslant \lambda_0 \leqslant L$, at $\mu$, if $\lambda_0 < \mu$ and $L$, if $\lambda_0 > L$. That is, the condition on imaginary part can be written as:
    \begin{equation*}
        |y| \leqslant
        \begin{cases}
            \sqrt{\alpha^2 \mu^2 - (x - \mu)^2}, \quad \frac{x}{1 - \alpha^2} < \mu, \\
            \frac{\alpha}{\sqrt{1 - \alpha^2}} x, \quad \mu \leqslant \frac{x}{1 - \alpha^2} \leqslant L, \\
            \sqrt{\alpha^2 L^2 - (L - x)^2}, \quad \frac{x}{1 - \alpha^2} > L.
        \end{cases}
    \end{equation*}
\end{proof}

\begin{lemma} \label{lem:polynomial_expression}
    Let an algorithm $\mathcal{A}$ satisfy~\eqref{eq:linear_first_order}, let $f$ be quadratic~\eqref{eq:quadratic_problem} with matrix $A$ and let $\wtgg f$ be given by the affine oracle~\eqref{eq:affine_quadratic_model} with matrix $\dA$ ($\wtA = A + \dA$). Then, for every $N \geqslant 0$, there exists a polynomial $p_N$ such that the iterate $x^N$ and $p_N$ satisfy:
    \begin{equation*}
        x^N - x^* = p_N(\wtA) (x^0 - x^*), \quad p_N(0) = 1, \quad \deg p_N \leqslant N.
    \end{equation*}
\end{lemma}
\begin{proof}
    Let us denote $y^k := x^k - x^*$. Since $f$ is a quadratic function with matrix $A$, its gradient is $\nabla f(x^k) = A(x^k - x^*) = A y^k$. By the condition of the affine oracle with matrix $\wtA = A + \dA$, we have $\wtgg f(x^k) = \wtA(x^k - x^*) = \wtA y^k.$ We will prove by induction. 
    
    Base case $N = 0$. We have
    \begin{equation*}
        y^0 = x^0 - x^* = I (x^0 - x^*) = p_0(\wtA)(x^0 - x^*),
    \end{equation*}
    where $p_0(\lambda) = 1$. Clearly, $\deg p_0 = 0$ and $p_0(0) = 1$.
    
    Inductive step. Assume that for all $k = 0, \dots, N$ it has been proven that $y^k = p_k(\wtA)y^0, \quad p_k(0) = 1, \quad \deg p_k \leqslant k$. From the linear first-order condition~\eqref{eq:linear_first_order}, there exist coefficients $h_0, \dots, h_N$ such that $y^{N+1} = y^0 + \sum\limits_{k=0}^N h_k \wtA \, y^k$. Substituting the induction hypothesis into this expression yields:
    \begin{equation*}
        y^{N+1}
        = y^0 + \sum_{k=0}^N h_k \wtA \, p_k(\wtA)y^0
        = \left(I + \wtA \sum_{k=0}^N h_k p_k(\wtA)\right)y^0.
    \end{equation*}
    Then $y^{N+1} = p_{N+1}(\wtA)y^0$, where $p_{N+1}(\lambda) = 1 + \lambda \sum\limits_{k=0}^N h_k p_k(\lambda)$. Moreover $\deg p_{N+1} \leqslant 1 + \max\limits_{0 \leqslant k \leqslant N} \deg p_k \leqslant N+1$ and $p_{N+1}(0) = 1$. This completes the induction.
\end{proof}

\section{Complex analysis} \label{appendix:complex_misc}

\subsection{Lens properties} \label{appendix:lens_prop}

\begin{definition} \label{def:lens}
For $0 < a < b$ and $0 < \theta < \frac{\pi}{2}$, we define a lens as a result of a conformal map:
\begin{equation}\label{eq:lens_def}
    \lens(a,b,\theta) = \left\{ a + (b - a) \frac{\omega}{1 + \omega} \st
     |\Arg(\omega)| \leqslant \theta \text{ or } \omega = \infty \text{ or } \omega = 0 \right\}.
\end{equation}
That is:
\begin{equation*}
    \begin{gathered}
        g(\omega) = \frac{a + b \omega}{1 + \omega}, \quad \mathcal{A} = \left\{ \omega \in \wC \st |\Arg(\omega)| \leqslant \theta \text{ or } \omega = \infty \text{ or } \omega = 0 \right\}, \\
        \lens(a, b, \theta) = g(\mathcal{A}).
    \end{gathered}
\end{equation*}
\end{definition}

\begin{lemma} \label{lem:lens_circl_intersect}
    Let $0 < a < b$ and $0 < \theta < \frac{\pi}{2}$ and $\lens(a, b, \theta)$ be the lens (Definition~\ref{def:lens}). Then:
    \begin{equation*}
        \begin{gathered}
            \lens(a, b, \theta) = D_1 \cap D_2, \quad D_1 = \{ z \in \C \st |z - c_1 | \leqslant r \}, \quad D_2 = \{ z \in \C \st |z - c_2 | \leqslant r \}, \\
            c_1 = \frac{a + b}{2} - i \frac{b - a}{2} \cot(\theta), \quad c_2 = \frac{a + b}{2} + i \frac{b - a}{2} \cot(\theta), \quad r = \frac{b - a}{2 \sin(\theta)}.
        \end{gathered}
    \end{equation*}
\end{lemma}
\begin{proof}
    Let us consider:
    \begin{equation*}
        \begin{gathered}
            g(\omega) = \frac{a + b \omega}{1 + \omega}, \quad
            \lens(a, b, \theta) = g(\mathcal{A}), \\
            \mathcal{A} = \left\{ \omega \in \wC \st |\Arg(\omega)| \leqslant \theta \text{ or } \omega = \infty \text{ or } \omega = 0 \right\}, \\
            A_1 = \{ \omega \in \wC \st \theta - \pi \leqslant \Arg(\omega) \leqslant \theta \text{ or } \omega = \infty \text{ or } \omega = 0 \}, \\
            A_2 = \{ \omega \in \wC \st -\theta \leqslant \Arg(\omega) \leqslant \pi - \theta \text{ or } \omega = \infty \text{ or } \omega = 0 \}, \\
            L_1 = \{ t e^{i \theta} \st t \in \R \} \cup \{ \infty \}, \quad L_2 = \{ t e^{-i \theta} \st t \in \R \} \cup \{ \infty \}
        \end{gathered}
    \end{equation*}
    With fractional linear mapping $g$, line $L_1$ will transform into a circle, since the pole of this map $\omega_0 = -1, \; \omega_0 \not\in L_1$. The center $c_1$ of this circle will be symmetric to $\infty$ with respect to the circle itself, which means that point $\omega_0$ is symmetric to point $g^{-1}(c_1)$ with respect to line $L_1$, that is $g^{-1}(c_1) = -e^{2i\theta} \Leftrightarrow c_1 = g(-e^{2i\theta})$.
    \begin{eqnarray*}
        c_1 & = & g(-e^{2i\theta}) = a + (b - a) \frac{-e^{2i\theta}}{1 - e^{2i\theta}} = a + (b - a) \frac{-e^{i\theta}}{e^{-i\theta} - e^{i\theta}} 
        \\
        & = & a + (b - a) \frac{-\cos \theta - i \sin \theta}{\cos \theta - i \sin \theta - \cos \theta - i \sin \theta}
        \\
        & = & a + (b - a) \left(\frac{1}{2} - \frac{i}{2} \cot \theta  \right) = \frac{a + b}{2} - i \frac{b - a}{2} \cot \theta.
    \end{eqnarray*}
    The radius of the corresponding circle can be computed, since $a$ belongs to it:
    \begin{equation*}
        r = \sqrt{\left(\frac{a + b}{2} - a\right)^2 + \frac{(b - a)^2}{4} \cot^2 \theta}
        = \frac{b - a}{2} \sqrt{1 + \cot^2 \theta} = \frac{b - a}{2 \sin \theta}.
    \end{equation*}
    Since the pole $\omega_0 \not\in A_1$, then $g(A_1) = D_1$. Similarly $g(A_2) = D_2$. Using injectivity of $g$, one can obtain:
    \begin{equation*}
        \lens(a, b, \theta) \overset{\text{Definition}~\ref{def:lens}}{=} g(\mathcal{A}) = g(A_1 \cap A_2) = g(A_1) \cap g(A_2) = D_1 \cap D_2.
    \end{equation*}
\end{proof}

\begin{lemma} \label{lem:admissible_lens}
    Let $0 < \mu \leqslant L, 0 < \alpha < 1$, and:
    \begin{equation*}
        \begin{gathered}
            \perm = \{ z \in \C \st \exists \lambda \in [\mu; L], \; |z - \lambda | \leqslant \alpha \lambda \}, \\
            a = (1 - \alpha) \mu, \quad b = (1 + \alpha)L, \quad \theta = \arcsin \alpha,
        \end{gathered}
    \end{equation*}
    then $\lens(a, b, \theta) \subseteq \perm$.
\end{lemma}
\begin{proof}
    From Lemma~\ref{lem:lens_circl_intersect}, one can obtain:
    \begin{equation*}
        \begin{gathered}
            \lens(a, b, \theta) = D_1 \cap D_2, \quad D_1 = \{ z \in \C \st |z - c_1 | \leqslant r \}, \quad D_2 = \{ z \in \C \st |z - c_2 | \leqslant r \}, \\
            c_1 = x_0 - i y_0, \quad c_2 = x_0 + i y_0, \quad r = \frac{b - a}{2 \sin(\theta)}, \quad x_0 = \frac{a + b}{2}, \quad y_0 = \frac{b - a}{2} \cot(\theta).
        \end{gathered}
    \end{equation*}
    Then $\lens(a, b, \theta)$ can be expressed as:
    \begin{equation*}
        \begin{gathered}
            \lens(a, b, \theta) = \{z \in \C \st |\Imm z | \leqslant s(\Ree z), \; a \leqslant \Ree z \leqslant b \}, \\
            s(x) = \sqrt{r^2 - \left(x - \frac{a + b}{2}\right)^2} - \frac{b - a}{2} \cot \theta.
        \end{gathered}
    \end{equation*}
    Using Lemma~\ref{lem:permissible_spectrum_envelope}:
    \begin{equation*}
        \begin{gathered}
            \perm = \{ z \in \C \st |\Imm z | \leqslant u(\Ree z), a \leqslant \Ree z \leqslant b \}, \\
            u(x) = 
            \begin{cases}
                \sqrt{\mu^2 \sin^2(\theta) - (x - \mu)^2}, \quad x < \mu \cos^2\theta, \\
                x \tan \theta, \quad \mu \cos^2 \theta \leqslant x \leqslant L \cos^2 \theta, \\
                \sqrt{L^2 \sin^2(\theta) - (L - x)^2}, \quad x > L \cos^2 \theta.
            \end{cases}
        \end{gathered}
    \end{equation*}
    Then, to prove $\lens(a, b, \theta) \subseteq \perm$, it is sufficient to show that $a \leqslant x \leqslant b, \; s(x) \leqslant u(x)$. Firstly:
    \begin{equation*}
        s'(x) = \displaystyle \frac{\frac{a + b}{2} - x}{\sqrt{r^2 - \left(x - \frac{a + b}{2}\right)^2}},
    \end{equation*}
    and $s(x)$ is concave (its graph is an upper circular arc of $\partial D_1$), then $s(x) \leqslant s(a) + s'(a) (x - a)$. Since $s(a) = 0, \; s'(a) = \tan \theta$, then $s(x) \leqslant (x - a) \tan \theta \leqslant x \tan \theta$. That is $s(x) \leqslant u(x)$ for all $\mu \cos^2 \theta \leqslant x \leqslant L \cos^2 \theta$. 

    Taking into account $a \leqslant x < \mu \cos^2\theta$, inequality $s(x) \leqslant u(x)$ is equivalent to studying the relative positions of $\partial D_1$ and $\{z \in \C \st |z - \mu| \leqslant \alpha \mu \}$ arcs. It is a standard fact in geometry that two distinct circles cannot have more than two points of intersection. Since $s(a) = u(a) = 0$, these two circles already have $1$ common point. Also $s(\mu \cos^2\theta) \leqslant \mu \cos^2\theta \tan \theta = u(\mu \cos^2\theta)$, $u(x) = \sqrt{(x - a)((1 + \sin \theta)\mu - x)}$:
    \begin{equation*}
        \lim_{x \to a+} \frac{s(x)}{u(x)} \leqslant \lim_{x \to a+} \frac{(x - a) \tan \theta}{\sqrt{(x - a)((1 + \sin \theta)\mu - x})} = 0,
    \end{equation*}
    then we have $s(x) < u(x)$ for all $x > a$ sufficiently close to $a$. Suppose that $s(\widehat{x}) > u(\widehat{x})$ for some $\widehat{x} \in (a, \mu\cos^2\theta)$. We can then choose $\widetilde{x}\in(a,\widehat{x})$ such that $s(\widetilde{x})<u(\widetilde{x})$. By continuity and the inequality at $\mu\cos^2\theta$,
    there exist $\widetilde{x}< x_1 <\widehat{x}< x_2 \le\mu\cos^2\theta$ such that $s(x_1)=u(x_1)$ and $s(x_2)=u(x_2)$. Together with the intersection at $a$, this gives three distinct intersection points of two distinct circles, a contradiction.

    Similarly, considering the case $L \cos^2\theta \leqslant x \leqslant b$, one can note that $s(b) = u(b) = 0$, $s(L \cos^2 \theta) \leqslant L \cos^2 (\theta) \tan(\theta) = u(L \cos^2 \theta)$, $s'(b) = -\tan \theta \Rightarrow s(x) \leqslant (b - x) \tan \theta$ and:
    \begin{equation*}
        \lim_{x \to b-} \frac{s(x)}{u(x)} \leqslant \lim_{x \to b-} \frac{(b - x) \tan \theta}{\sqrt{(b - x)(x - (1 - \sin \theta)L)}} = 0.
    \end{equation*}
    Then $s(x) \leqslant u(x)$ for $L \cos^2\theta \leqslant x \leqslant b$. Since $(\forall x \in [a; b]) \; s(x) \leqslant u(x)$, we can conclude $\lens(a, b, \theta) \subseteq \perm$.
\end{proof}

\begin{lemma} \label{lem:rieman_map_lens}
    Let $0 < a < b$ and $0 < \theta < \frac{\pi}{2}$ and $\lens(a, b, \theta)$ be the lens (Definition~\ref{def:lens}). Then there is a conformal map $\Phi$, such that:
    \begin{equation*}
        \begin{gathered}
            \Bc = \{ z \in \C \st |z| > 1 \} \cup \{ \infty \}, \quad \Phi(z) = \frac{\left(\frac{a - z}{b - z} \right)^{\gamma} + 1}{\left(\frac{a - z}{b - z} \right)^{\gamma} - 1}, \\
            \gamma = \frac{\pi}{2 ( \pi - \theta)}, \quad \Phi\left( \wC \setminus \lens(a, b, \theta) \right) = \Bc.
        \end{gathered}
    \end{equation*}
\end{lemma}
\begin{proof}
    From the definition of the lens (Definition~\ref{def:lens}) it is the image of a sector on the extended complex plane. That is:
    \begin{equation*}
        \begin{gathered}
            f_1(z) = \frac{a - z}{b - z}, \quad \mathcal{A}_1 = \left\{ \omega \in \wC \st |\Arg(\omega)| \geqslant \pi - \theta \text{ or } \omega = \infty \text{ or } \omega = 0 \right\}, \\
            \lens(a, b, \theta) = f_1^{-1}(\mathcal{A}_1), \quad f_1^{-1}(\omega) = \frac{b \omega - a}{\omega - 1}, \\
            \wC \setminus \mathcal{A}_1 = f_1(\wC \setminus \lens(a, b, \theta)) \text{ , since it is a linear fractional transformation}.
        \end{gathered}
    \end{equation*}
    Set $\wC \setminus \mathcal{A}_1 = \left\{ \omega \in \wC \st |\Arg(\omega)| < \pi - \theta \right\}$ can be transformed (choosing principal value of the corresponding function) to the right half-plane:
    \begin{equation*}
        \begin{gathered}
            f_2(z) = z^{\gamma} = \exp(\gamma \Log z), \quad \gamma = \frac{\pi}{2( \pi - \theta)}, \\
            \mathcal{A}_2 = \left\{ \omega \in \wC \st |\Arg(\omega)| < \frac{\pi}{2} \right\}, \quad \mathcal{A}_2 = f_2(\wC \setminus \mathcal{A}_1).
        \end{gathered}
    \end{equation*}
    The standard transformation of the half-plane $\mathcal{A}_2$ into a set $\Bc$ is:
    \begin{equation*}
        \begin{gathered}
            f_3(z) = \frac{z + 1}{z - 1}, \quad \Bc = f_3(\mathcal{A}_2).
        \end{gathered}
    \end{equation*}
    Then we can build the target map:
    \begin{equation*}
        \Phi(z) = (f_3 \circ f_2 \circ f_1) (z) = \frac{\left(\frac{a - z}{b - z} \right)^{\gamma} + 1}{\left(\frac{a - z}{b - z} \right)^{\gamma} - 1}, \quad \gamma = \frac{\pi}{2 ( \pi - \theta)}, \quad \Phi(\wC \setminus \lens(a, b, \theta)) = \Bc.
    \end{equation*}
    Since each of the functions $f_1, f_2, f_3$ is conformal on its domain, the map $\Phi$ is also conformal.
\end{proof}

\subsection{Potential theory} \label{appendix:complex_potential}

\begin{definition}[Definition 3.1.1~\cite{ransford1995potential}] \label{def:potential}
    Let $\mu$ be a finite Borel measure on $\C$ with compact support. Its potential is the function $p_\mu \colon \C \to [-\infty, \infty)$ defined by
    \begin{equation*}
        p_\mu(z) = \int_{\C} \log |z - w| \, d\mu(w).
    \end{equation*}
\end{definition}

\begin{definition}[Definition 3.2.1~\cite{ransford1995potential}] \label{def:energy}
    Let $\mu$ be a finite Borel measure on $\C$ with compact support. Its energy $I(\mu)$ is given by
    \begin{equation*}
        I(\mu) = \int_{\C} p_\mu(z) \, d\mu(z).
    \end{equation*}
\end{definition}

\begin{definition}[Definition 3.2.2~\cite{ransford1995potential}] \label{def:polar_set}
    A subset $E$ of $\C$ is called polar if $I(\mu) = -\infty$ for every finite Borel measure $\mu \neq 0$ for which $\supp \mu$ is a compact subset of $E$.
\end{definition}

\begin{proposition}[Corollary 3.2.4~\cite{ransford1995potential}] \label{prop:polar_borel_set}
    Every Borel polar set has Lebesgue measure zero.
\end{proposition}

\begin{proposition}[Corollary 3.6.4~\cite{ransford1995potential}] \label{prop:polar_disconnected_set}
    Every closed polar set E is totally disconnected.
\end{proposition}

\begin{definition} \label{def:green_function}
    Let $D$ be a proper subdomain of $\wC$. A Green's function for $D$ is a map $g_D \colon D \times D \to (-\infty, +\infty]$, such that for each $w \in D$:
    \begin{enumerate}
        \item[(a)] \(g_D(\cdot, w)\) is harmonic on \(D \setminus \{w\}\), and bounded outside each neighbourhood of \(w\);
        \item[(b)] \(g_D(w, w) = \infty\), and as \(z \to w\),
        \[
            g_D(z, w) = \begin{cases}
                \log|z| + O(1), & w = \infty, \\
                -\log|z - w| + O(1), & w \neq \infty;
            \end{cases}
        \]
        \item[(c)] $g_D(z, w) \to 0$ as $z \to \zeta$, for all $\zeta \in \partial D$, except some polar set.
    \end{enumerate}
\end{definition}

\begin{proposition} \label{prop:exterior_ball_green_infty}
    Let $\Bc = \{z \st |z| > 1 \} \cup \{ \infty \}$, then the Green's function is given by:
    \begin{equation*}
        g_{\Bc}(z, \infty) = \ln |z|,
    \end{equation*}
\end{proposition}
\begin{proof}
    All properties of Green's function can be checked directly.
\end{proof}

\begin{theorem}[Theorem 4.4.4~\cite{ransford1995potential}] \label{thm:green_function_map}
   Let $D_1$ and $D_2$ be domains in $\wC$ with non-polar boundaries, and let $f: D_1 \to D_2$ be a conformal map onto $D_2$. Then:
   \begin{equation*}
        g_{D_2}(f(z), f(w)) = g_{D_1}(z, w) \quad (z, w \in D_1).
   \end{equation*}
\end{theorem}

\begin{lemma} \label{lem:green_function_using_conf_map}
    Let $D$ be a proper subdomain with non-polar boundary of $\wC$, such that $\infty \in D$. $\Phi: D \to \Bc$ (where $\Bc = \{z \st |z| > 1 \} \cup \{ \infty \}$) is a conformal map between $D$ and $\Bc$, such that $\Phi(\infty) = \infty$. Then the Green's function is given by:
    \begin{equation*}
        g_D(z, \infty) = \ln|\Phi(z)|,
    \end{equation*}
\end{lemma}
\begin{proof}
    Consequence of Theorem~\ref{thm:green_function_map} and Proposition~\ref{prop:exterior_ball_green_infty}.
\end{proof}

\section{Lower bound} \label{appendix:lower_bound}

The following finite extremal construction is closely related to classical max--min and min--max approximation results for normal matrices (see~\cite{liesen2013max}). We include a self-contained proof in the form needed below.

\begin{theorem} \label{thm:existance_points_z_k}
    Let $K \subset \C$ be an infinite compact set. Let $N \geqslant 0$ and
    \begin{equation*}
         \mathcal P_N = \{P \in \R[z]: \deg P \leqslant N,\ P(0) = 1\},
         \quad
         M = \inf_{P \in \mathcal P_N} \max_{z\in K}|P(z)|.
    \end{equation*}
    Then there exist nodes $z_1,\dots,z_{N + 1}\in K$, and $\omega_1,\dots,\omega_{N + 1} \geqslant 0$ such that $\sum\limits_{k=1}^{N + 1} \omega_k = 1$ and:
    \begin{equation}
        (\forall P \in \mathcal P_N) \quad \sum_{k=1}^{N + 1} \omega_k |P(z_k)|^2 \geqslant M^2.
    \end{equation}
\end{theorem}
\begin{proof}
    Let us denote $F(a, z) = \left|1 + \sum_{k = 1}^{N} a_k z^k \right|^2$ -- the squared modulus of the polynomial $P \in \P$  with coefficients $a \in \R^N$ at point $z \in K$. That is $M^2 = \inf\limits_{a \in \R^N} \max\limits_{z \in K} F(a, z)$. One can note that for any fixed $z \in K$ the function $F(a, z)$ is a quadratic function of the variable $a$. Then we can consider the function:
    \begin{equation*}
        \varphi(a) = \max_{z \in K} F(a, z).
    \end{equation*}
    Let us introduce a function:
    \begin{equation*}
        s(a) = \left\| \sum_{k = 1}^{N} a_k z^k \right\|_K = \max_{z \in K} \left | \sum_{k = 1}^{N} a_k z^k \right |.
    \end{equation*}
    One can see that:
    \begin{equation*}
        \varphi(a) = \max_{z \in K} \left|1 + \sum_{k = 1}^{N} a_k z^k \right|^2
        = \left \|1 + \sum_{k = 1}^{N} a_k z^k \right\|_K^2
        \geqslant \left( \left \| \sum_{k = 1}^{N} a_k z^k \right \|_K - 1 \right)^2 = \left(s(a) - 1 \right)^2.
    \end{equation*}
    It should be noted that $s(a) = 0 \Leftrightarrow a = 0$ (since $K$ is infinite), so $s$ is a norm on $\R^N$. That is:
    \begin{equation*}
        \lim_{\|a \|_2 \to +\infty} s(a) = +\infty \Rightarrow \lim_{\|a \|_2 \to +\infty} \varphi(a) = +\infty,
    \end{equation*}
    that is $\varphi$ is continuous and coercive. Since $F(\cdot, z)$ is convex and differentiable one can apply Theorem~\ref{thm:danskin_sub}:
    \begin{equation*}
        \partial \varphi(a) = \conv \{\nabla_a F(a, z) \mid z \in Z(a)\}, \quad Z(a) = \left\{ \overline{z} \;\middle|\; F(a, \overline{z}) = \max_{z \in K} F(a, z) \right\}.
    \end{equation*}
    Since $\varphi$ is continuous and coercive, $M^2 = \inf_a\varphi(a)$ is attained at some $a^*\in\mathbb R^N$. Since $a^*$ minimizes $\varphi$, Theorem~\ref{thm:subdiff_optim} can be applied: $0 \in \partial \varphi(a^*)$. From Theorem~\ref{thm:caratheodory_thm} applied to the set $\{\nabla_a F(a^*, z) \mid z \in Z(a^*)\}$: $\exists z_1, \dots z_{N + 1} \in Z(a^*), \exists \omega_1, \dots \omega_{N + 1}$:
    \begin{equation*}
        \sum_{k = 1}^{N + 1} \omega_k \nabla_a F(a^*, z_k) = 0, \quad \sum_{k = 1}^{N + 1} \omega_k = 1, \quad \omega_k \geqslant 0, \quad F(a^*, z_k) = \varphi(a^*) = M^2.
    \end{equation*}
    Then we can consider the function:
    \begin{equation*}
        G(a) = \sum_{k = 1}^{N + 1} \omega_k F(a, z_k).
    \end{equation*}
    Since $G$ is a quadratic convex function and $\nabla G(a^*) = 0$, then $(\forall P \in \P)$:
    \begin{equation*}
        M^2 = \sum_{k = 1}^{N + 1} \omega_k F(a^*, z_k) = G(a^*) \leqslant G(a) = \sum_{k = 1}^{N + 1} \omega_k F(a, z_k) = \sum_{k = 1}^{N + 1} \omega_k | P(z_k) |^2.
    \end{equation*}
\end{proof}

\begin{theorem} \label{thm:lens_polynom_max}
    Let $0 < a < b$ and $0 < \theta < \frac{\pi}{2}$ and let $\lens(a, b, \theta)$ be the lens (Definition~\ref{def:lens}). Then:
    \begin{equation*}
        \begin{gathered}    
            \| P(z) \|_{\lens(a, b, \theta)} \geqslant \left( \frac{1 - (\nicefrac{a}{b})^{\gamma}}{1 + (\nicefrac{a}{b})^{\gamma}}\right)^N, \quad \gamma = \frac{\pi}{2 ( \pi - \theta)}, \\
            P \in \P = \{P \in \R[z]: \deg P \leqslant N,\ P(0) = 1 \}.
        \end{gathered}
    \end{equation*}
\end{theorem}
\begin{proof}
    Since $\lens(a, b, \theta)$ has positive Lebesgue measure, then it is not polar (Proposition~\ref{prop:polar_borel_set}). $\partial \lens(a, b, \theta)$ is connected (it is the union of two arcs of two circles by Lemma~\ref{lem:lens_circl_intersect}), then it is also not polar by Proposition~\ref{prop:polar_disconnected_set}. Then using Lemma~\ref{lem:rieman_map_lens} and Lemma~\ref{lem:green_function_using_conf_map} we can obtain Green's function for $D = \wC \setminus \lens(a, b, \theta)$:
    \begin{equation*}
        \begin{gathered}
            \Phi(z) = \frac{\left(\frac{a - z}{b - z} \right)^{\gamma} + 1}{\left(\frac{a - z}{b - z} \right)^{\gamma} - 1}, \quad
            \gamma = \frac{\pi}{2 ( \pi - \theta)}, \quad
            g_{D}(z, \infty) = \ln |\Phi(z) |.
        \end{gathered}
    \end{equation*}
    If $\deg P = 0 \Longrightarrow P \equiv 1$, that is the desired inequality is immediate. If $\deg P \geqslant 1$, then Theorem~\ref{thm:bernstein_theorem} can be applied for $P \in \P$:
    \begin{equation*}
        (\forall z \in D \setminus \{\infty\}) \quad \left( \frac{|P(z)|}{\|P\|_{\lens(a, b, \theta)}} \right)^{1/\deg P} \leqslant e^{g_D(z, \infty)} = | \Phi(z) |.
    \end{equation*}
    Substituting $z = 0$:
    \begin{equation*}
        \frac{1}{\| P \|^{1/\deg P}_{\lens(a, b, \theta)}} \leqslant | \Phi(0) | = \frac{1 + (\nicefrac{a}{b})^{\gamma}}{1 - (\nicefrac{a}{b})^{\gamma}},
    \end{equation*}
    then:
    \begin{equation*}
         \| P \|_{\lens(a, b, \theta)} \geqslant \left( \frac{1 - (\nicefrac{a}{b})^{\gamma}}{1 + (\nicefrac{a}{b})^{\gamma}}\right)^{\deg P} \geqslant \left( \frac{1 - (\nicefrac{a}{b})^{\gamma}}{1 + (\nicefrac{a}{b})^{\gamma}}\right)^{N}.
    \end{equation*}
\end{proof}

\begin{proof}[Proof of Theorem~\ref{thm:quadratic_lower_bound_single_point}]
    Let us consider $K = \lens((1 - \alpha) \mu, (1 + \alpha)L, \arcsin \alpha)$ and:
    \begin{equation*}
        \begin{gathered}
             \P = \{P \in \R[z]: \deg P \leqslant N,\ P(0) = 1\},
             \quad
             M = \inf_{P \in \P} \max_{z\in K}|P(z)|.
        \end{gathered}
    \end{equation*}
    Since $K \subset \C$ is compact we can apply Theorem~\ref{thm:existance_points_z_k}, that is $\exists z_1 \dots z_{N + 1} \in K, \; \omega_1,\dots,\omega_{N + 1} \geqslant 0$ such that:
    \begin{equation*}
        (\forall P \in \P) \quad \sum_{k=1}^{N + 1} \omega_k |P(z_k)|^2 \geqslant M^2, \quad \sum_{k = 1}^{N + 1} \omega_k = 1.
    \end{equation*}
    Using Lemma~\ref{lem:admissible_lens} -- $K \subseteq \perm$, where $\perm = \{ z \in \C \st \exists \lambda \in [\mu; L], \; |z - \lambda | \leqslant \alpha \lambda \}$, that is for all $1 \leqslant k \leqslant N + 1, \; \exists \lambda_k \in [\mu; L], \; |\lambda_k - z_k| \leqslant \alpha \lambda_k$. Let us consider the homomorphism:
    \begin{equation*}
        T: \C \to \R^{2 \times 2}, \quad
        T(z) = \begin{pmatrix} \Ree z & - \Imm z \\ \Imm z & \Ree z \end{pmatrix},
    \end{equation*}
    and matrices:
    \begin{equation*}
        A = \diag(T(\lambda_1), \dots T(\lambda_{N + 1})), \quad \wtA = \diag(T(z_1), \dots T(z_{N + 1})).
    \end{equation*}
    Lemma~\ref{lem:lambda_and_z} states that $\dA = \wtA - A$ satisfies the relative error condition~\eqref{eq:affine_quadratic_model}. Let $f(x) = \frac{1}{2} x^T A x$; then $f^* = 0$ and $x^* = 0$. We choose $x^0 = (x^0_{1}, \dots, x^0_{N + 1}) \in \R^{2(N + 1)}$, where each $x^0_k \in \R^2$ is defined as $x^0_k = \sqrt{\frac{\omega_k}{2}} R (1, 1)$. Lemma~\ref{lem:polynomial_expression} states that any first-order method~\eqref{eq:linear_first_order} $\mathcal{A}$ for the function $f$ can be expressed as $x^N = P(\wtA) x^0$, where $P \in \P$. Then:
    \begin{eqnarray*}
        \|x^N - x^* \|_2^2 
        & = & \| P(\wtA) x^0 \|_2^2 = \sum_{k = 1}^{N + 1} \|T(P(z_k)) x^0_k \|_2^2 \overset{\text{Lemma}~\ref{lem:T_operator_isometry}}{=} \sum_{k = 1}^{N + 1} | P(z_k) |^2 \| x^0_k \|_2^2 
        \\
        & = & R^2 \sum_{k = 1}^{N + 1} \omega_k | P(z_k) |^2 \geqslant R^2 M^2.
    \end{eqnarray*}
    Using Theorem~\ref{thm:lens_polynom_max} we obtain the lower bound:
    \begin{equation*}
        \begin{gathered}    
        \max_{z \in K} |P(z)| \geqslant \left( \frac{1 - (\nicefrac{a}{b})^{\gamma}}{1 + (\nicefrac{a}{b})^{\gamma}}\right)^N, \\
        \gamma = \frac{\pi}{2 ( \pi - \theta)}, \quad a = (1 - \alpha) \mu, \quad b = (1 + \alpha) L, \quad \theta = \arcsin \alpha.
        \end{gathered}
    \end{equation*}
    Then:
    \begin{eqnarray*}
        \|x^N - x^* \|_2^2 
        & \geqslant & R^2 M^2 = R^2 \inf_{P \in \P} \max_{z\in K}|P(z)|^2 \geqslant R^2 \left( \frac{1 - (\nicefrac{a}{b})^{\gamma}}{1 + (\nicefrac{a}{b})^{\gamma}}\right)^{2 N}
        \\
        & = & R^2 \left(\frac{b^{\gamma} - a^{\gamma}}{b^{\gamma} + a^{\gamma}} \right)^{2 N}
    \end{eqnarray*}
\end{proof}

\end{appendix}

\end{document}